\documentclass[reqno,11pt]{amsart}

\usepackage{amsthm, amsmath, amssymb, stmaryrd}
\usepackage[utf8]{inputenc}
\usepackage[T1]{fontenc}

\usepackage[hyphenbreaks]{breakurl}
\usepackage[hyphens]{url}
\usepackage{systeme}
\usepackage[shortlabels]{enumitem}
\usepackage[hidelinks]{hyperref}
\usepackage{microtype}

\usepackage{bm}
\usepackage[margin=1in]{geometry}

\usepackage[textsize=scriptsize,backgroundcolor=orange!5]{todonotes}

\usepackage[noabbrev,capitalize,sort,nameinlink]{cleveref}
\crefname{equation}{}{}
\numberwithin{equation}{section}

\usepackage{mathtools}
\newtheorem{theorem}{Theorem}[section]
\newtheorem{proposition}[theorem]{Proposition}
\newtheorem{lemma}[theorem]{Lemma}
\newtheorem{claim}[theorem]{Claim}

	\crefname{claim}{Claim}{Claims}

\theoremstyle{definition}
\newtheorem{definition}[theorem]{Definition}

\theoremstyle{remark}

\newcommand{\abs}[1]{\left\lvert#1\right\rvert}

\newcommand{\ceil}[1]{\left\lceil #1 \right\rceil}

\newcommand{\set}[1]{\left\{ #1 \right\}}

\newcommand*{\PP}{\mathbb{P}}

\newcommand{\cF}{\mathcal F}

\newcommand{\veps}{\varepsilon}

\newlength{\hght}

\makeatletter
\newcommand\thankssymb[1]{\textsuperscript{\@fnsymbol{#1}}}
\makeatother

\author[Ting-Wei Chao]{Ting-Wei Chao\thankssymb{1}}
\author[Bruce Reed]{Bruce Reed\thankssymb{2}}

\thanks{\thankssymb{1}Department of Mathematics, Massachusetts Institute of Technology, Cambridge, MA 02139\@. Email: {\tt twchao@mit.edu}}

\thanks{\thankssymb{2}Mathematical Institute, Academia Sinica, Taiwan. Supported by  NSTC Grant 115-2115-M-001 -010\@. Email: {\tt bruce.al.reed@gmail.com }}

\title{A nearly tight bound on the probability $G(n,p)$ is even-degenerate}

\begin{document}

\maketitle
\begin{abstract}
    A graph is even-degenerate if we can order its vertices as $v_1,\dots,v_n$ so that the number of edges from $v_i$ to $v_1,\dots,v_{i-1}$ is even for all $i>2$. We show that the probability $G(n,p)$ is even-degenerate is $e^{-\Omega(n/\log\log n)}$, which is only a factor of $\log\log n$ away, in the exponent, from the trivial lower bound $2^{-O(n)}$.
\end{abstract}
\section{Introduction}
A graph $G$ is {\it even-degenerate} if its vertices can be ordered as $(v_1,\dots,v_n)$ so that every vertex except $v_2$ has an even number of neighbors appearing earlier in the order. The study of such graphs was initiated by Janzer and Yip \cite{JY24}, as they are a special case of even-decomposable graphs whose study was initiated by Alon \cite{Alon24} as an outgrowth of his work on graph codes, see also \cite{AlonGKMS23,Ver25}. We note the exception for $v_2$ is to permit graphs with an odd number of edges to be even-degenerate. 

The focus of work on even-degenerate graphs has been to bound the probability that the random graph $G(n,p)$ is even-degenerate.
$G(n,p)$ is the graph on vertex set  $\{1,\dots,n\}$ constructed by taking $\binom{n}{2}$ independent identically distributed random $\{0,1\}$ variables, one for each pair of vertices, which are $1$ with probability $p$, as indicator variables for the existence of the edges. 

As pointed out by Janzer and Yip, if $n$ is even, the probability that all vertices of $G(n,p)$ have odd degree is at least $2^{-O(n)}$ while if $n$ is odd the probability that one vertex has degree $0$ and the rest have odd degree is also at least $2^{-O(n)}$. So, the probability that $G(n,p)$ is not even-degenerate is at least $2^{-O(n)}$. Exploiting symmetry arising from the fact that for $p=\frac{1}{2}$ we have $p=1-p$ they were able to show that for $p=\frac{1}{2}$ the probability that $G(n,p)$ is not even-degenerate is at most $2^{-\Omega(n)}$. They asked if this was true for all $p$. Chao, Dong, and Xu \cite{CDX26} showed that for every $\alpha>0$, for any constant $p$ with $0<p<1$ the probability that $G(n,p)$ is not even-degenerate is at most $2^{-\Omega_\alpha(n^{1/2-\alpha})}$.

We improve this result by showing:
\begin{theorem}\label{thm:main}
    For any constant $p\in (0,1)$, we have
    \[\PP(G(n,p)\textup{ is not even-degenerate})=e^{-\Omega_p(n/\log\log n)}.\]
\end{theorem}

In our proof, we build the sequence $(v_1,\dots,v_n)$ recursively, similarly to Janzer--Yip and Chao--Dong--Xu. The main difference is that we build the sequence from $v_1$ to $v_n$, whereas the previous approaches build it from $v_n$ to $v_1$. The advantage is that, when we look for the next vertex $v_i$, we only need to query each candidate vertex $u$, which exposes only the edges between $u$ and $v_1,\dots,v_{i-1}$. This saves a substantial amount of randomness. However, since $p$ may be far from $\frac{1}{2}$, a queried vertex can be added with probability as small as $q=\min(p,1-p)$. This leads us to a multi-round scanning process, which we call the \emph{greedy scanning algorithm}. In each step of the recursion, we can apply the greedy scanning algorithm with a set $B_1$ of vertices set aside, which gives a recursion of the form $1-f(n)\geq(1-f(m))(1-\veps)$, where $f$ is the failure probability. We can bootstrap this idea by setting aside several disjoint sets $B_1,\dots,B_b$ instead of one. By showing that the processes for different $B_i$ are conditionally independent, we obtain a recursion of the form $f(n)\leq f(m)^b+\veps$. The proof is completed by proving the base case and solving the recursions.

\textbf{Notation.} Throughout this paper, $p\in(0,1)$ is fixed, and $q=\min(p,1-p)$.
For a vertex $v$ and a set $S$ of vertices, write $d_S(v)=|N(v)\cap S|$. Let $e(H)$ be the number of edges of a graph $H$.

\textbf{AI usage.} The main proof ideas are human generated. ChatGPT 6 pro was used for a preliminary calculation to figure out the optimal parameters. Claude Opus 5.5 max and ChatGPT 6 pro were also used to simplify the arguments and polish the writing. In particular, the technical definition, \cref{def:instance}, was suggested by ChatGPT to make the recursion cleaner.
\section{Upper bound}

We apply an iterative recursion
where in each iteration we consider extending an already constructed prefix of a proposed order to a longer prefix. This requires us to  consider a generalization of even degeneracy, where 
we specify for each vertex the desired parity of its degree into the vertices appearing earlier in the order.
For technical simplicity for any input graph $G$ we first choose a prefix $v_1,v_2$ of length 2. Any such prefix can be extended 
into  an ordering certifying the even degenerateness of $G$ if and only if there is an ordering $w_1,\dots,w_{|V(G)|-2}$ of $G-\{v_1,v_2\}$
such that for every vertex $w_j$, $d_{\{v_1,v_2\}}(w_j)=d_{\{w_1,\dots,w_{j-1}\}}(w_j) \pmod{2}$. To extend this idea, suppose that we have already fixed a prefix, we will consider the random graph on the remaining vertices, and labels that carry the parities from the prefix.

A \emph{labeled graph} is a pair $(H,\lambda)$ consisting of a graph $H$ and a map $\lambda\colon V(H)\to\{0,1\}$. If $\lambda(v)=0$, we may \emph{select} $v$. That is, we delete $v$ and flip the labels $\lambda(u)$ of its neighbors $u\in N(v)$, so that the label becomes $1-\lambda(u)$. We say that $(H,\lambda)$ is \emph{good} if there is an ordering of $V(H)=\{v_1,v_2,\dots, v_{\abs{V(H)}}\}$ so that we can sequentially select them, from $v_1$ to $v_{\abs{V(H)}}$, and we call the sequence $(v_1,\dots,v_{\abs{V(H)}})$ a \emph{good sequence} of $(H,\lambda)$.
Note that $(H,\lambda)$ can only be good if $e(H)+\sum_v \lambda(v)=0\pmod{2}$. Thus, we define the \emph{parity} of $(H,\lambda)$ to be $e(H)+\sum_{v}\lambda(v)\pmod{2}$. We know that the parity is invariant under selecting vertices. We note that $G$ has an order certifying even-degeneracy beginning $v_1,v_2$ if and only if the labeled graph with $H=G-\{v_1,v_2\}$ and $\lambda(v)=d_{\{v_1,v_2\}}(v)]\pmod{2}$ is good. 

We will consider an intermediate step where the vertices from some $X \subseteq V(H)$  have been selected, except for a leftover set $U$ of vertices whose labels are all $1$, and $B$ is a set of vertices the edges to which have yet to be examined. By applying a \emph{greedy scanning algorithm}, we will reduce to a smaller intermediate step where there are only a few choices of the leftover set $U$ for the next round. The key idea is that we can run this greedy scanning algorithm several times, each time preserving a different set $B_i$ as the set of unused vertices. We will pick a family of possible $U$ for the next round that does not depend on which round we are running, so the failure events of the next round will be conditionally independent. To carry out this idea, we begin with the following definition.

\begin{definition}\label{def:instance}
Let $n\ge1$. An \emph{instance} $(G,X,\cF,B,(\xi_w)_{w\in B})$ of size $n$ consists of
\begin{enumerate}[(i)]
\item a random graph $G$ on the vertex set $V(G)=X\sqcup B$,
\item the set $X$ of at most $n^2$ vertices with $G[X]$ being deterministic,
\item a family $\cF\subseteq 2^X$ with $\abs{\cF}\leq e^{n^{3/4}}$,
\item the set $B$ of $n$ vertices disjoint from $X$ and independent Bernoulli random variables, which are the indicators of edges in $G$, associated with $X \times B$ and the unordered pairs of $B$ each of which is one with probability $p$,
\item a random parity $\xi_w$ for each $w\in B$, which we call \emph{bits}, that are mutually independent from each other and from the variables of (iv) such that each is 1 with some probability between $p$ and $1-p$.
\end{enumerate}
We remark that each $U\in \cF$ is a possible leftover from the previous round, and for each $w\in B$, the bit $\xi_w$ represents the parity of the number of edges from $w$ to all the vertices that are already selected outside of $V(G)$. Since $X\setminus U$ is also selected and $U$ can not be selected in the previous round, we define 
\begin{align*}
    \lambda_U(w)=\begin{cases}
1&\textup{ if $w\in U$}\\
\xi_w+d_{X\setminus U}(w) \pmod{2}&\textup{ if $w\in B$}
\end{cases}
\end{align*}
to be the current label as if $U$ is indeed the leftover. Let $\Gamma_U=(G[U\cup B],\lambda_U)$ be the labeled graph.
The instance \emph{fails} if there exists $U\in\cF$ such that the parity of $\Gamma_U$ is even but $\Gamma_U$ is not good. Otherwise, we say that the instance \emph{succeeds}. Let $f(n)$ be the maximum of the failure probability over all instances of size $n$.
\end{definition}

In this definition, one should think of $X$ as the set of vertices used in the last round, $U$ as the set of leftovers, and $\lambda$ as the parity of the number of edges to the vertices outside of $U\cup B$. 

We will show that $f(n)= e^{-\Omega(n/\log\log n)}$ by showing the following recursion and base case.
\begin{lemma}\label{lem:recursion}
There is a constant $s_0$ such that the following holds. Let $s,m,b$ be positive integers with $s_0\leq s\leq m^{4/3}$, and let $n=s+bm$ and $\veps=e^{2n^{3/4}}(1-q)^{s/2}$. Then
\begin{enumerate}[1.]
    \item $f(n)\leq f(m)^{b}+\veps$, and
    \item $1-f(n)\geq\left(1-f(m)\right)(1-\veps)$.
\end{enumerate}
\end{lemma}

\begin{lemma}\label{lem:start}
    There is a constant $n_0$ such that $f(n)<1$ for every $n\geq n_0$.
\end{lemma}

We will first prove these two lemmas, and then show that they together imply \cref{thm:main}.

\subsection{Greedy scanning algorithm}
We first state our greedy scanning algorithm, which will be crucial in the proofs of \cref{lem:recursion,lem:start}.

Let $W$ and $A$ be disjoint sets of vertices with $\abs{A}=s$ and $\abs{W\cup A}=N$, and fix an ordering of $W\cup A$ in which every vertex of $W$ comes before every vertex of $A$. Let $(H,\lambda)$ be a random labeled graph on $W\cup A$ such that
\begin{itemize}
\item the graph $H[W]$ and the labels $\lambda(v)$ for $v\in W$ are fixed;
\item the pairs of vertices with at least one end in $A$ are edges independently with probability $p$;
\item the labels $\lambda(a)$ for $a\in A$ are independent Bernoulli random variables with probability between $p$ and $1-p$, and are independent of the edges.
\end{itemize}

The \emph{greedy scanning algorithm} on $(H,\lambda)$ runs as follows.
\begin{itemize}
    \item In the $i$-th round with $1\leq i\leq N+1$, we will find a remaining vertex $v_i$ with label $\lambda(v_i)=0$ and select it as follows.
    \begin{itemize}
        \item We run through the remaining vertices in the prescribed ordering of $W\cup A$ one by one. When we move to a vertex $u$, we expose its label $\lambda(u)$. We call this a \emph{query} of $u$.
        \item If $\lambda(u)=0$, then we take $v_i=u$. Otherwise, we move to the next vertex.
        \item If there is no such vertex or if all the vertices are selected, then the algorithm ends.
    \end{itemize}
    \item Assume the algorithm ends in the $k$-th round. Then it outputs a sequence of selected vertices $(v_1,\dots,v_{k-1})$ and a set of remaining vertices $R_0=(W\cup A)\setminus \{v_1,\dots,v_{k-1}\}$, where all the labels in $R_0$ become $1$ after selecting $(v_1,\dots,v_{k-1})$.
\end{itemize}
Note that to expose $\lambda(u)$, we expose only the edges between $u$ and $v_1,\dots,v_{i-1}$, and the initial randomness of $\lambda(u)$ if $u$ is queried for the first time.

We will show that $R_0$ is likely to be a subset of $A$ of a special form. For an integer $t\geq 0$, let 
\[\cF_t(A)=\left\{R\subseteq A\,\middle|\, \abs{R\Delta S}\leq 2t\textup{ for some suffix $S$ of $A$}\right\}.\] 
We first bound the size of $\cF_t(A)$.
\begin{lemma}\label{lemma:sizeFtA}
    We have $\abs{\cF_t(A)}\leq (s+1)^{2t+1}.$
\end{lemma}
\begin{proof}
    Each $R\in\cF_t(A)$ is determined by a suffix $S$ of $A$ and the set $R\Delta S$, which has at most $2t$ vertices. There are $s+1$ possible suffixes of $A$, and at most $(s+1)^{2t}$ subsets of $A$ with at most $2t$ vertices. Thus, $\abs{\cF_t(A)}\leq (s+1)^{2t+1}.$
\end{proof}

Now, we show that with high probability $R_0\in\cF_t(A)$.
\begin{lemma}\label{lem:scan}
There is a constant $s_1$ such that the following holds. If $s\geq s_1$ and $t=\ceil{\sqrt{s}}$, then the set $R_0$ of remaining vertices output by the greedy scanning algorithm on $(H,\lambda)$ satisfies
\[\PP(R_0\notin\cF_t(A))\leq N^2(1-q)^{s/2}.\]
\end{lemma}

\begin{proof}
If $R_0\neq\varnothing$, then the algorithm ends in a round in which no remaining vertex has label $0$. In this final round, every vertex of $R_0$ is queried and has label $1$. Moreover, if $u$ comes before $v_i$ in the ordering and $u\notin\set{v_1,\dots,v_{i-1}}$, then $u$ is queried in the $i$-th round.

Note that when a vertex $u\in W\cup A$ is queried in the $i$-th and the $j$-th round with $i<j$, but not any other rounds in between, then the edge $uv_i$ is still random in the $j$-th round query, unless both $u$ and $v_i$ are in $W$. Also, if $u$ is queried for the first time and $u\in A$, then the initial label is exposed for the first time. In both cases, the probability that $\lambda(u)=0$ conditional on the process until this query is between $p$ and $1-p$.

\begin{claim}\label{claim:W}
    For each $x\in W$,
    \[\PP(x\in R_0)\leq N(1-q)^s.\]
\end{claim}
\begin{proof}
    For each $k\in [N]$, let $E_k$ be the event that the algorithm ends in the $k$-th round and $x\in R_0$. If $x\in R_0$, then at most $N-1$ vertices are selected, so the algorithm ends in the $k$-th round for some $k\in[N]$. It is sufficient to show that $\PP(E_k)\leq (1-q)^s$ for each $k$, and the claim follows from the union bound.

    Write $A\setminus R_0=\{u_1,\dots,u_r\}$ where $u_j$ is selected in the $i_j$-th round with $i_1<i_2<\dots<i_r<k$. We write $i_{r+1}=k$ for convenience. If $x\in R_0$, then $x$ must be queried but not selected in the $i_j$-th round for all $j\in [r]$. For each $j\in [r]$, take $i'\in (i_j,i_{j+1}]$ to be the smallest one such that $x$ is queried in the $i'$-th round. In this round, the edge $xu_j$ is still random, so $x$ is not selected with probability at most $1-q$. 
    
    The event that the algorithm ends in the $k$-th round tells us that none of the remaining vertices in $A\cap R_0$ is selected. Note that every $y\in A\cap R_0$ must be queried in the $k$-th round, and it is not selected with probability at most $1-q$. 

    In the discussion above, we considered $s$ queries:
    \begin{itemize}
        \item for each $u_j$, the next query of $x$ after $u_j$ is selected, and
        \item every query of a vertex of $A$ in the $k$-th round.
    \end{itemize}
    Each query gives a factor of at most $1-q$ to the event $E_k$, so $\PP(E_k)\leq (1-q)^s$.
\end{proof}

Next, we show that, with high probability, there is no large selected set coming after a large unselected set in the prescribed ordering.
\begin{claim}\label{claim:Q}
    Let $Q_1,Q_2\subseteq A$ be two sets of at least $t+1$ vertices such that every vertex of $Q_2$ comes after every vertex of $Q_1$ in the ordering. Then
    \[\PP(Q_1\subseteq R_0, Q_2\cap R_0=\varnothing)\leq (1-q)^{(t+1)^2}.\]
\end{claim}
\begin{proof}
    We may assume without loss of generality that $\abs{Q_1}=\abs{Q_2}=t+1$. Let $E$ be the event that $Q_1\subseteq R_0$ and $Q_2\cap R_0=\varnothing$. Suppose that $E$ occurs. Write $Q_2=\{u_1,\dots,u_{t+1}\}$ where $u_j$ is selected in the $i_j$-th round with $i_1<i_2<\dots<i_{t+1}$, let $k$ be the round in which the algorithm ends, and write $i_{t+2}=k$ for convenience. Let $x\in Q_1$. Since $x$ comes before every vertex of $Q_2$ and $x\in R_0$, $x$ must be queried but not selected in the $i_j$-th round for all $j\in [t+1]$, and also in the $k$-th round. For each $j\in [t+1]$, take $i'\in (i_j,i_{j+1}]$ to be the smallest one such that $x$ is queried in the $i'$-th round. In this round, the edge $xu_j$ is still random, so $x$ is not selected with probability at most $1-q$.

    In the discussion above, we considered $(t+1)^2$ queries, namely, for each $x\in Q_1$ and each $u_j$, the next query of $x$ after $u_j$ is selected. Each query gives a factor of at most $1-q$ to the event $E$, so $\PP(E)\leq (1-q)^{(t+1)^2}$.
\end{proof}

Now, we can finish the proof. Suppose that $R_0\notin\cF_t(A)$. If $R_0\not\subseteq A$, then $x\in R_0$ for some $x\in W$. Otherwise, we have $\abs{R_0}=\abs{R_0\Delta\varnothing}>2t$ since $\varnothing$ is a suffix of $A$. Let $Q_1$ be the set of the first $t+1$ vertices of $R_0$, and let $S$ be the suffix of $A$ whose first vertex is the last vertex of $Q_1$. Then $R_0\setminus S$ consists of $t$ vertices, so $\abs{S\setminus R_0}>t$. Let $Q_2$ be a set of $t+1$ vertices in $S\setminus R_0$. Every vertex of $Q_2$ comes after every vertex of $Q_1$, and we know that $Q_1\subseteq R_0,Q_2\cap R_0=\varnothing$. Thus, either there is a vertex $x\in W\cap R_0$, or there are two sets $Q_1,Q_2$ satisfying the assumption in \cref{claim:Q} with $Q_1\subseteq R_0, Q_2\cap R_0=\varnothing$. There are $N-s$ choices of $x$ and at most $s^{2t+2}$ choices of $(Q_1,Q_2)$. Thus, by \cref{claim:W,claim:Q} and the union bound, we have
\[\PP(R_0\notin\cF_t(A))\leq (N-s)N(1-q)^s+s^{2t+2}(1-q)^{(t+1)^2}.\]
When $s$ is sufficiently large, we have $(t+1)^2\geq s$ and $s^{2t+2}\leq (1-q)^{-s/2}$, since $(2t+2)\log s=O(\sqrt{s}\log s)$. Thus, the right-hand side is at most $\left((N-s)N+1\right)(1-q)^{s/2}\leq N^2(1-q)^{s/2}$.
\end{proof}

\subsection{The recursion}
In this subsection, we prove the recursion, \cref{lem:recursion}.
\begin{proof}[Proof of \cref{lem:recursion}]
    Fix an instance $(G,X,\cF,B,\xi_w)$ of size $n$, and partition $B$ into sets $A,B_1,\dots,B_b$ with $\abs{A}=s$ and $\abs{B_1}=\dots=\abs{B_b}=m$. Set $t=\ceil{\sqrt{s}}$. We fix an arbitrary ordering of each of the sets $X$, $A$, and $B_1,\dots,B_b$. For each $1\leq i\leq b$, we will run the greedy scanning algorithm on all the vertices except $B_i$, and then use $B_i$ as the set of unused vertices of an instance of size $m$.
    
    For $U\in\cF$ and $i\in [b]$, let $\Gamma_{U,i}$ be the labeled graph obtained from $\Gamma_U$ by deleting $B_i$ and keeping the labels of the other vertices. Note that a sequence of selections in $\Gamma_{U,i}$ is also a sequence of selections in $\Gamma_U$, which will also flip some labels in $B_i$.

    We run the greedy scanning algorithm on $\Gamma_{U,i}$, and let $D_{U,i}$ be the event that the set of remaining vertices is not in $\cF_t(A)$. Let $W_i=U\cup(\bigcup_{j\neq i}B_j)$. We may apply \cref{lem:scan} with $W_i$ in place of $W$, the set $A$, and the label $\lambda_U$. To fix the randomness in $W_i$, we need to condition on the edges inside $W_i$, the labels $\lambda_U(v)$ for $v\in W_i$, and the edges between $A$ and $X\setminus U$ first. If $s\geq s_1$, then the conditional probability of $D_{U,i}$ is at most $(n^2+n)^2(1-q)^{s/2}$ since the total number of vertices is at most $n^2+n$. In particular, we have $\PP(D_{U,i})\leq (n^2+n)^2(1-q)^{s/2}$. Let $D_i$ be the event that $D_{U,i}$ happens for some $U\in\cF$. Then by the union bound, we know that $\PP(D_i)\leq e^{n^{3/4}}(n^2+n)^2(1-q)^{s/2}$, and hence 
    \[\sum_{i=1}^b\PP(D_i)\leq be^{n^{3/4}}(n^2+n)^2(1-q)^{s/2}\leq e^{2n^{3/4}}(1-q)^{s/2}=\veps\]
    when $n$ is large enough, which is true if $s$ is large enough.

    Now, fix $i\in[b]$ and $U\in\cF$, and suppose that $D_i$ does not happen. Then the algorithm on $\Gamma_{U,i}$ selects every vertex of $W_i$, and leaves a set $R\in\cF_t(A)$ of vertices of label $1$. Since the label of a vertex $w\in B_i$ is flipped once for each selected neighbor of $w$, its current label is     $\xi'_w+d_{A\setminus R}(w)\pmod{2}$, where $\xi'_w=\xi_w+d_X(w)+\sum_{j\neq i}d_{B_j}(w)\pmod{2}$, which does not depend on $U$. So the remaining labeled graph is the labeled graph $\Gamma_R$ of $(G[A\cup B_i],A,\cF_t(A),B_i,\xi'_w)$, defined as in \cref{def:instance}.

     Let $Z$ be the collection of 
    \begin{itemize}
        \item the edges inside $A$,
        \item the bits $(\xi_a)_{a\in A}$, 
        \item the edges between $X$ and $B$, and 
        \item the edges between $B_j$ and $B_{j'}$ for all $j\neq j'$.
    \end{itemize} 
    Conditioned on $Z$, $(G[A\cup B_i],A,\cF_t(A),B_i,\xi'_w)$ is an instance of size $m$. Indeed, $\abs{A}=s\leq m^2$, it is not hard to check that the random variables satisfy the definition of instances, and the family $\cF_t(A)$ has at most $(s+1)^{2t+1}= e^{O(m^{2/3}\log m)}\leq e^{m^{3/4}}$ sets by \cref{lemma:sizeFtA} when $s$ is large enough.  Define $E_i$ to be the event that $(G[A\cup B_i],A,\cF_t(A),B_i,\xi'_w)$ fails. Then $\PP(E_i\mid Z)\leq f(m)$ by the definition of $f(m)$. It is clear that if neither $D_i$ nor $E_i$ happens, then $(G,X,\cF,B,\xi_w)$ succeeds. Indeed, if $\Gamma_U$ has even parity for some $U\in\cF$, then the remaining labeled graph $\Gamma_R$ also has even parity since the parity is invariant under selecting vertices, so we can combine the sequence of selecting vertices in the greedy scanning step and a good sequence of $\Gamma_R$.

    We note that the event $E_i$ depends only on the edges inside $A\cup B_i$, the edges between $B_i$ and $X\cup\bigcup_{j\neq i}B_j$, and the bits $(\xi_w)_{w\in B_i}$.

    Thus, given $Z$, the event $E_i$ depends only on the edges and bits within $B_i$, and edges from $B_i$ to $A$. So, conditioned on $Z$, the events $E_1,\dots,E_b$ are independent.
    
    Now, we prove part 1 of the lemma. If $(G,X,\cF,B,\xi_w)$ fails, then for every $i\in[b]$, either $D_i$ or $E_i$ happens. So either all of $E_1,\dots,E_b$ happen, or $D_i$ happens for some $i$. Since $E_1,\dots,E_b$ are conditionally independent given $Z$, we have
    \[\PP(E_1\cap\dots\cap E_b\mid Z)=\prod_{i=1}^b\PP(E_i\mid Z)\leq f(m)^b.\]
    By taking expectation over $Z$ and the union bound, we have
    \[\PP(\textup{$(G,X,\cF,B,\xi_w)$ fails})\leq f(m)^b+\sum_{i=1}^b\PP(D_i)\leq f(m)^b+\veps.\]
    This proves part 1. 

    For part 2, we consider conditioning on 
    $D_1$ not occurring and note that this event is the union of 
    some of the choices for edges and bits within $X\cup (\bigcup_{j\neq 1}B_j)\cup A$. Summing the conditional 
    expectation of $E_1$ over every such choice,
    we see that the conditional chance of $E_1$ not occurring 
    given $D_1$ does not occur is at least $1-f(m)$. 
    This proves part 2.
\end{proof}

\subsection{The base case}
In this subsection, we prove the base case, \cref{lem:start}.
\begin{proof}[Proof of \cref{lem:start}]
    Let $s_1$ be the constant in \cref{lem:scan}, let $n$ be sufficiently large, and fix an instance $(G,X,\cF,B,\xi_w)$ of size $n$. Let $B=A\sqcup B_0$ with $\abs{A}=\lfloor \frac{n-1}{2} \rfloor$. We will run the greedy scanning algorithm with $B_0$ set aside, and use the edges incident  $B_0$ to complete every possible set of remaining vertices.

    For $U\in\cF$, let $\Gamma'_U$ be the labeled graph obtained from $\Gamma_U$ by deleting $B_0$ and keeping the labels of the other vertices, and let $R_U$ be the set of remaining vertices output by the greedy scanning algorithm on $\Gamma'_U$ with $W=U$ and the set $A$ in a fixed ordering. As in the proof of \cref{lem:recursion}, after conditioning on the edges between $A$ and $X\setminus U$, \cref{lem:scan} applies with $s=\abs{A}$ and $N\leq n^2+n$. Since $\abs{A}\geq\frac{n}{3}\geq s_1$ and every member of $\cF_t(A)$ is a subset of $A$, where $t=\ceil{\sqrt{s}}$, the union bound over $\cF$ gives
\[\PP(R_U\not\subseteq A\textup{ for some }U\in\cF)\leq e^{n^{3/4}}(n^2+n)^2(1-q)^{s/2}<1,\]
    where the last inequality holds since $s=\abs{A}\geq\frac{n}{3}$ and $n$ is large. Let $E_1$ be the event that $R_U\subseteq A$ for every $U\in\cF$. Then $\PP(E_1)>0$, and $E_1$ is determined by the graph on $X\cup A$ and the bits $(\xi_w)_{w\in A}$.

    Now, it is sufficient to show the following claim.
    \begin{claim}\label{claim:base}
        There exists a deterministic assignment of the edges within $B_0$ and from $B_0$ to $X\cup A$, and the bits $\xi_w$ for $w\in B_0$, such that the following holds for any assignment of edges in $A$. For any $R\subseteq A$, let $\lambda'_R$ be the label with $\lambda'_R(w)=1$ for $w\in R$ and $\lambda'_R(w)=\xi_w+d_{A\setminus R}(w)\pmod{2}$ for $w\in B_0$. Then the labeled graph $(G[R\cup B_0],\lambda'_R)$ is good if it has even parity.
    \end{claim}
     Once this is shown, we know that the event $E_2$ that this assignment happens has positive probability. Moreover, suppose that $E_1$ and $E_2$ happen, and $\Gamma_U$ has even parity for some $U\in\cF$. After selecting $U\cup(A\setminus R_U)$ by the greedy scanning algorithm, the label of each $w\in B_0$ becomes $\xi_w+d_{X}(w)+d_{A\setminus R_U}(w)$, so the remaining labeled graph is $(G[R_U\cup B_0],\lambda'_R)$ in the claim with $R=R_U$, which also has even parity. So we can combine the sequence of selecting $U\cup(A\setminus R_U)$ given by $E_1$ and the sequence of selecting $R_U\cup B_0$ given by $E_2$. This gives a good sequence of $\Gamma_U$. Since the random edges and bits that $E_1$ and $E_2$ depend on are disjoint, we have
    \[\PP((G,X,\cF,B,\xi_w)\textup{ succeeds})\geq \PP(E_1\cap E_2)=\PP(E_1)\PP(E_2)>0.\]
    Since the instance is arbitrary, and the maximum in the definition of $f$ is achievable since the space is compact, this implies that $f(n)<1$ for every sufficiently large $n$, which proves \cref{lem:start}.

    It remains to prove the claim. 
    \begin{proof}[Proof of \cref{claim:base}.]
        We can write  $A=\set{x_1,\dots,x_r}$ and the ordered set $B_0$ as $\set{y_1,\dots,y_r,z}$ followed possibly by a final vertex $z^*$.  We assign the edges and bits as follows:
\begin{enumerate}[(i)]
\item the edge within $B_0$ form a star with center $z$
and leaves $\{y_1,\dots,y_r\}$.  

\item the edges between $X \cup A$ and $B_0$ form a matching $\{x_1y_1,\dots,x_ry_r\}$;
\item $\xi_w=0$ for every $w\in B_0-z$.
\end{enumerate}
    Fix any assignments of edges inside $A$ resulting from performing our greedy scanning algorithm with $B_0$ removed for the given choice as discussed above, and fix a choice of $R \subseteq A$.  
    
    The first time we consider selecting $y_i$, we have not yet selected $z$, so by (i), (ii), and (iii) we select $y_i$ precisely if we have not yet selected $x_i$. But selecting $y_i$ changes the labels of only $z$ and $x_i$ so the next vertex selected will be $x_i$. 
    Thus before considering selecting $z$ for the first time, 
    we have selected all of $A$. 

    Now when we first consider selecting $z$, by (iii) $z^*$ has label $0$, and a vertex other than $z$ or $z^*$  has label $0$ precisely if it is selected. Since $z$ sees all the unselected vertices except $z^*$, if 
    the labeled graph has even parity we select $z$. Since the
    unselected vertices form a  subset of the stable set
    $B_0-z$, (iii) implies we now select each of them in turn. 
    
    \end{proof}
\end{proof}

\subsection{Proof of \cref{thm:main}}

Let $n_0$ be the constant given by \cref{lem:start}. We may assume that $n_0\geq 16$, so that $\log\log n\geq 1$ for $n\ge n_0$. All the constants below depend only on $p$. We first show that $f(n)$ is bounded away from $1$.
\begin{proposition}
    There is a constant $c_0>0$ such that $f(n)\leq 1-c_0$ for every $n\geq n_0$.
\end{proposition}
\begin{proof}
We apply part 2 of \cref{lem:recursion} with $b=1$ repeatedly. Let $s=\ceil{n^{4/5}}$ and $g(n)=n-\ceil{n^{4/5}}$. When $n$ is sufficiently large, we have $s_0\leq s\leq g(n)^{4/3}$, $g(n)\geq\frac{n}{2}$, and $e^{2n^{3/4}}(1-q)^{s/2}\leq\frac{1}{n^2}$. So part 2 of \cref{lem:recursion} gives
    \[1-f(n)\geq\left(1-\frac{1}{n^2}\right)\left(1-f(g(n))\right).\]
    Take $N_0\geq 2n_0$ such that this holds for every $n\geq N_0$.

    Let $n\geq N_0$. We apply the inequality above repeatedly, until we reach a number $m'=g^{(k)}(n)<N_0$ for some $k>0$, where $g^{(k)}$ is the $k$-fold composition of $g$. Since $g(x)\geq \frac{x}{2}$ for every $x\geq N_0$, we have $m'\geq\frac{N_0}{2}\geq n_0$. The numbers $n,g(n),\dots,g^{(k-1)}(n)$ are all distinct and at least $2$, so we have
    \[1-f(n)\geq\prod_{j=2}^{\infty}\left(1-\frac{1}{j^2}\right)\cdot\left(1-f(m')\right)=\frac{1}{2}\left(1-f(m')\right).\]
    Thus, we may pick $c_0=\frac{1}{2}\min_{m'\in [n_0,N_0]}(1-f(m'))$, which is positive because of \cref{lem:start}.
\end{proof}

Next, we prove the upper bound on $f(n)$ we need.
\begin{proposition}\label{prop:f}
There is a constant $c>0$ such that $f(n)\leq e^{-cn/\log\log n}$ for every $n\geq n_0$.
\end{proposition}

\begin{proof}
We will apply part 1 of \cref{lem:recursion} to $f(n)$ by choosing the parameters as follows. For $n\geq 16$, let $m=\ceil{n^{3/4}}$, let $b$ be the integer with
\[
  bm\leq n-\frac{n}{16\log\log n}<(b+1)m,
\]
and let $s=n-bm$. Then
\[\frac{n}{16\log\log n}\leq s<\frac{n}{16\log\log n}+m.\]
Note that when $n$ is sufficiently large, we have $s_0\leq s\leq n\leq m^{4/3}$, $m\geq n_0$, and
\[e^{2n^{3/4}}(1-q)^{s/2}\leq e^{-2an/\log\log n}\]
for some $a>0$. We also have $m\leq 2n^{3/4}\leq n^{4/5}$, and hence $\log\log m\leq \log\log n-\log\frac{5}{4}$. Since $(b-1)m>n-\frac{n}{16\log\log n}-2m$, we have
\[\frac{(b-1)m}{\log\log m}\geq \frac{n-\frac{n}{16\log\log n}-2m}{\log\log n-\log\frac{5}{4}}\geq \frac{n}{\log\log n}\]
because $\frac{n}{\log\log n}>(\frac{n}{16\log\log n}+2m)/\log\frac{5}{4}$ when $n$ is large enough. Thus, 
\[\frac{bm}{\log\log m}\geq \frac{n}{\log\log n}+\frac{m}{\log\log m}\]
in this case.
We will also need the fact that
\[e^{-(c_0/n)\cdot m/\log\log m}+e^{-an/\log\log n}\leq 1\]
when $n$ is large enough. This is because 
\[1-e^{-(c_0/n)\cdot m/\log\log m}\geq \frac{1}{2}\cdot\frac{c_0m}{n\log\log m},\]
which decreases slower than $e^{-an/\log\log n}$ as $n\rightarrow\infty$.

Now, we choose $n_1\geq n_0$ such that for every $n\geq n_1$ the conclusions above hold, and also $c_0/n_1\leq a$. Set $c=c_0/n_1$. It follows that
\[
  f(n)\leq 1-c_0\leq e^{-c_0}\leq e^{-cn_1}\leq e^{-cn/\log\log n}
\]
for all $n_0\leq n\leq n_1$.
Now let $n>n_1$, and assume the bound for all numbers in $[n_0,n)$. Since $n_0\leq m<n$, the induction hypothesis and the choice of $n_1$ give
\[
  f(m)^{b}\leq e^{-cbm/\log\log m}\leq e^{-cn/\log\log n-cm/\log\log m}.
\]
By part 1 of \cref{lem:recursion}, $c\leq a$, and $c\geq c_0/n$, we have
\[
  f(n)\leq f(m)^{b}+e^{-2an/\log\log n}\leq e^{-cn/\log\log n}\left(e^{-cm/\log\log m}+e^{-an/\log\log n}\right)\leq e^{-cn/\log\log n},
\]
and we are done by induction.
\end{proof}
Now, we are ready to prove \cref{thm:main} using \cref{prop:f}.
\begin{proof}[Proof of \cref{thm:main}]
    
Recall that we want to find an ordering $v_1,\dots,v_n$ of the vertices in $V(G)$ such that for every $i>2$ the vertex $v_i$ has an even number of neighbors among $v_1,\dots,v_{i-1}$.

Consider the random graph $G=G(n,p)$ with $n\geq n_0+2$, and fix a choice of $v_1$. We expose all the edges incident to $v_1$, and fix two vertices $v_{2,0},v_{2,1}$ so that $v_1v_{2,0}\notin E(G)$ and $v_1v_{2,1}\in E(G)$ if they exist. Both vertices exist unless $v_1$ is isolated or adjacent to all other vertices, which happens with probability $p^{n-1}+(1-p)^{n-1}\leq 2(1-q)^{n-1}$.

For $j\in \{0,1\}$, let $B_j=V(G)\setminus\set{v_1,v_{2,j}}$, and give each $w\in B_j$ the bit $\xi_{w,j}=d_{\set{v_1,v_{2,j}}}(w) \pmod2$.
Since we exposed all the edges incident to $v_1$, $d_{\set{v_1}}(w)$ is fixed and the edge $wv_{2,j}$ is present independently with probability $p$. So the bits are independent Bernoulli random variables with probability between $p$ and $1-p$, and are independent of the edges of $G[B_j]$. Thus, $B_j$ and $(\xi_{w,j})_{w\in B_j}$ form an instance of size $n-2$ with $X=\varnothing,\cF=\set{\varnothing}$, so the label $\lambda_\varnothing(w)=\xi_{w,j}$ for all $w\in B_j$. If $(w_1,\dots,w_{n-2})$ is a good sequence of $(G[B_j],\lambda_\varnothing)$, then the sequence $v_1,v_{2,j},w_1,\dots,w_{n-2}$ shows that $G$ is even-degenerate.

If both $v_{2,0}$ and $v_{2,1}$ exist, then for $j\equiv e(G)\pmod2$ the labeled graph $(G[B_j],\lambda_\varnothing)$ has even parity. Hence, if $G$ is not even-degenerate, then either some $v_{2,j}$ does not exist, or for some $j$ the vertex $v_{2,j}$ exists and the corresponding instance fails. By \cref{prop:f},
\[
  \PP(G\text{ is not even-degenerate})\leq 2(1-q)^{n-1}+2f(n-2)\leq 2(1-q)^{n-1}+2e^{-c(n-2)/\log\log (n-2)},
\]
which is $e^{-\Omega(n/\log\log n)}$.
\end{proof}

\bibliographystyle{amsplain0}
\bibliography{bib}
\end{document}